%Paper: funct-an/9401002
%From: GUIDO@MAT.UTOVRM.IT
%Date: Tue, 18 JAN 94 18:21 GMT

%%%%%% FORMAT %%%%%%%%%%%%
\magnification\magstep1
\hoffset=0.5truecm
\voffset=0.5truecm
\hsize=15.8truecm
\vsize=23.truecm
\baselineskip=14pt plus0.1pt minus0.1pt \parindent=19pt
\lineskip=4pt\lineskiplimit=0.1pt      \parskip=0.1pt plus1pt

%%%%%% Extra Fonts %%%%%%%%%

 \def\got#1{\bf#1}

%% If eufm10 fonts are available, you will get a much more
%% gothic effect for Lie algebras by replacing previous line
%% with the following:
%% \font\got=eufm10 scaled\magstep1

\def\G{\hbox{\got g}}
\def\P{\hbox{\got p}}
\def\V{\hbox{\got a}}

%%%%%% SYMBOLS %%%%%%%%%%%%%
\def\A{{\cal A}}

\def\a{\alpha}

\def\b{\beta}

\def\Co{{\bf C}}
\def\d{\delta}
\def\D{\Delta}

\def\f{\varphi}

\def\H{{\cal H}}
\def\K{{\cal K}}
\def\L{\Lambda}

\def\Na{{\bf N}}
\def\o{\omega}
\def\O{{\cal O}}
\def\p{\pi}
\def\Q{\Omega}

\def\R{{\cal R}}
\def\Re{{\bf R}}
\def\s{\sigma}

\def\th{\theta}
\def\U{{\cal U}}
\def\W{{\cal W}}

\def\Ze{{\bf Z}}
\def\ad{{\rm ad}}

\def\Poi{{\cal P}_+^\uparrow}
\def\Spin{\widetilde{\cal P}}
\def\imply{\Rightarrow}
\def\np{\par\noindent}
\def\nu{\char'43}
\def\lar{\longrightarrow}
\def\ov{\overline}

\def\quot#1#2{\hbox{$\matrix{#1\cr\vphantom{#2}}$
    \hskip-.8cm\vbox{\vfill/\vfill}\hskip-12.4cm
    $\matrix{\vphantom{#1}\cr#2}$}}

%%%%%%%% Macros for theorems %%%%%%%%
\def\proof{\medskip\noindent{\bf Proof.}\quad}
\def\proofof#1{\medskip\noindent{\bf Proof of #1.}\quad}
\def\square{\hbox{$\sqcap\!\!\!\!\sqcup$}}
\def\endproof{\hskip-4mm\hfill\vbox{\vskip3.5mm\square\vskip-3.5mm}\bigskip}
\def\section #1\par{\vskip0pt plus.3\vsize\penalty-75
    \vskip0pt plus -.3\vsize\bigskip\bigskip
    \noindent{\sectionfont #1}\nobreak\smallskip\noindent}

\def\claim#1#2\par{\vskip.1in\medbreak\noindent{\bf #1.} {\sl #2}\par
    \ifdim\lastskip<\medskipamount\removelastskip\penalty55\medskip\fi}
\def\rmclaim#1#2\par{\vskip.1in\medbreak\noindent{\bf #1.} { #2}\par
    \ifdim\lastskip<\medskipamount\removelastskip\penalty55\medskip\fi}

%%%%%%%% AUTO REF. %%%%%%%%%%%%%
\font\sectionfont=cmbx10 scaled\magstep1
\newcount\REFcount \REFcount=1
\def\numref{\number\REFcount}
\def\addref{\global\advance\REFcount by 1}
\def\wdef#1#2{\expandafter\xdef\csname#1\endcsname{#2}}
\def\wdch#1#2#3{\ifundef{#1#2}\wdef{#1#2}{#3}
    \else\write16{!!doubly defined#1,#2}\fi}
\def\wval#1{\csname#1\endcsname}
\def\ifundef#1{\expandafter\ifx\csname#1\endcsname\relax}
%%%%%%%%%%%%%%%%%%%%%%%%%%%%%%%%%%%%%%%%%%%%
\def\ref(#1){\wdef{q#1}{yes}\ifundef{r#1}$\diamondsuit$#1
  \write16{!!ref #1 was never defined!!}\else\wval{r#1}\fi}
\def\inputreferences{
    \def\REF(##1)##2\endREF{\wdch{r}{##1}{\numref}\addref}
    \REFERENCES}
\def\references{
    \def\REF(##1)##2\endREF{
        \ifundef{q##1}\write16{!!ref. [##1] was never quoted!!}\fi
        \item{[\ref(##1)]}##2}
    \section{References}\par\REFERENCES}
%%%%%%REFERENCES%%%%%%%%
 \def\REFERENCES{
 \REF(BiWi1)Bisognano J., Wichmann E. , ``{\it On the duality
condition for a Hermitian scalar field}", J. Math. Phys. {\bf
16} (1975), 985-1007.
 \endREF
 \REF(Borc1)Borchers H.J. , ``{\it The CPT theorem in
two-dimensional theories of local observables}", Commun. Math.
Phys. {\bf 143} (1992), 315.
 \endREF
 \REF(Brow1)Brown K.S., ``{\it Cohomology of groups}" Springer,
New York Heidelberg Berlin, 1982.
 \endREF
 \REF(BGL)Brunetti R., Guido D., Longo R. , ``{\it Modular
structure and duality in conformal quantum field theory}",
Commun. Math. Phys., {\bf 156} (1993), 201-219.
 \endREF
 \REF(BuSu1)Buchholz, D. Summers, S.J. ``{\it An algebraic
characterization of vacuum states in Minkowski space}", preprint
DESY 92-119
 \endREF
 \REF(DoLo1)Doplicher S., Longo R., ``{\it Standard and split
inclusions of von Neumann algebras}", 	Invent. Math. {\bf 73}
(1984), 493-536.
 \endREF
 \REF(FGL)In preparation.
 \endREF
 \REF(FrGa1)Gabbiani F., Fr\"ohlich J., ``{\it Operator algebras
and Conformal Field Theory}", pre-print ETH Z\"urich, 1992.
 \endREF
 \REF(GL1)Guido, D. Longo, R., ``{\it Relativistic invariance
and charge conjugation in quantum field theory}", Commun. Math.
Phys. {\bf 148} 521-551.
 \endREF
  \REF(Haag1)Haag R., {\it Local Quantum Physics}, Springer
Verlag, Berlin Heidelberg 1992.
 \endREF
 \REF(HiLo1)Hislop P., Longo R., ``{\it Modular structure of the
von~Neumann algebras	associated with the free massless scalar
field theory}", 		Commun. Math. Phys. {\bf 84} (1982), 84.
 \endREF
 \REF(McLa1)Mc Lane S. ``{\it Homology}", Springer, Berlin, 1963.
 \endREF
 \REF(Mich1) Michel L., ``{\it Sur les extensions centrales du
groupe  de Lorentz inhomog\`ene connexe}'', Nucl. Phys. 57 (1965)
356-385.
 \endREF
 \REF(Miln1)Milnor J., ``{\it Introduction to algebraic
K-theory}", Ann. of Math Studies {\bf 72}, Princeton University
Press, Princeton, 1971.
 \endREF
 \REF(Moor1)Moore C.C., ``{\it Group extensions and cohomology
for locally compact groups. III}'', Trans. Amer. Math. Soc. {\bf
221} (1976) 1.
 \endREF
 \REF(Moor2)Moore C.C., ``{\it Group extensions and cohomology
for locally compact groups. IV}'', Trans. Amer. Math. Soc. {\bf
221} (1976) 35.
 \endREF
 \REF(Simm)  Simms D.J., ``{\it Lie groups and quantum
mechanics}" Lecture Notes in Math., Vol. 52, Springer, Berlin,
1968.
 \endREF
 \REF(StraZs1)Str\u atil\u a S., Zsido L., {\it Lectures on von
Neumann algebras}, Abacus press, England 1979.
 \endREF
 \REF(Wies1)Wiesbrock H.V., ``{\it A comment on a recent work of
Borchers}", to appear in Lett. Math. Phys.
 \endREF
 \REF(Wies2)Wiesbrock H.V., ``{\it Conformal quantum field
theory and half-sided modular inclusions of von~Neumann
algebras}", to appear in Commun. Math. Phys.
 \endREF
 \REF(Yngv1)Yngvason J. ``{\it A note on essential
duality}", preprint.
 % available via anonymous ftp in mp_arc@math.utexas.edu
 \endREF
 }

\inputreferences

\font\ftitle=cmbx10 scaled\magstep1

\vskip2truecm
\centerline{\ftitle  GROUP COHOMOLOGY, }
\medskip\centerline{\ftitle MODULAR THEORY}
\medskip\centerline{\ftitle AND SPACE-TIME SYMMETRIES}
 \bigskip\bigskip
 \centerline{R. Brunetti$^1$\footnote{$^*$}{ Supported in
part by MURST and
CNR-GNAFA.}
 \footnote{$^\bullet$}{Supported in part by INFN, sez.
Napoli and by Fondazione ``A. della Riccia"},
 D. Guido$^{2*}$ and R. Longo$^{23*}$}
\footnote{}{E-mail\ brunetti@na.infn.it,
guido@mat.utovrm.it, longo@mat.utovrm.it
}
 \vskip1.truecm
\item{$(^1)$}  Department of Physics, Syracuse University
 \par
201 Physics Building, Syracuse, NY 13244-1130.
 \par
 \item{$(^2)$} Dipartimento di Matematica, Universit\`a
di Roma ``Tor Vergata''
\par
Via della Ricerca Scientifica, I--00133 Roma,
Italia.
\item{$(^3)$} Accademia Nazionale dei Lincei,
\par via della Lungara 10, I--00165 Roma, Italia
 \vskip 3cm\noindent
{\bf Abstract. } The Bisognano-Wichmann property on
the geometric behavior of the modular group of the von Neumann
algebras of local observables associated to wedge regions in
Quantum Field Theory is shown to provide an intrinsic sufficient
criterion for the existence of a covariant action of the
(universal covering of) the Poincar\'e group. In particular
this gives, together with our previous results,  an intrinsic
characterization of positive-energy conformal pre-cosheaves of
von Neumann algebras. To this end we adapt to our use Moore
theory of central extensions of locally compact  groups by
polish groups, selecting and making an analysis of a wider class
of extensions with natural measurable properties and showing
henceforth that the universal covering of the Poincar\'e group
has only trivial central extensions (vanishing of the first and
second order cohomology) within our class.

\section{Introduction}
\par
 In this paper we discuss the problem of characterizing
the existence of a  Poincar\'e covariant action for a net of
local observable algebras in terms of net-intrinsic algebraic
properties.
 In other words, given a state with a Reeh-Schlieder cyclicity
property, we look for a condition that ensures it to be a
relativistic vacuum.

The work done by Bisognano and Wichmann [\ref(BiWi1)]
during the  seventies about the essential duality has shown the
geometrical character of the Tomita-Takesaki modular operators
associated with von Neumann algebras generated by Wightman
fields with support based on wedge regions and the vacuum
vector.

While an algebraic version of the Bisognano-Wichmann theorem in
the Poin\-ca\-r\'e covariant case is not yet available, there
are general results in this direction given by the recent work
of Borchers [\ref(Borc1)], see also
[\ref(BuSu1),\ref(Wies1),\ref(Yngv1)] for related results. Based
on the result of Borchers, there is a complete algebraic
analysis in [\ref(BGL)] and [\ref(FrGa1)] showing that, for a
conformally invariant theory, the modular groups associated with
 the von~Neumann algebras of double cones are associated with
special one-parameter subgroups of the conformal group leaving
invariant the given double cone (cf. [\ref(HiLo1)]).
Applications of the Bisognano-Wichmann theorem to the analysis
of the relations between the space-time symmetries and
particle-antiparticle symmetry in quantum field theory are
contained in [\ref(GL1)].

 Here we show how to reconstruct the Poincar\'e group
representation from the net-intrinsic property of {\it
modular covariance} (assumption 2.1), i.e. the modular operators
associated with the algebras of the wedge regions act
geometrically as boosts.

In our previous paper [\ref(BGL)] we showed in particular that,
for a Poincar\'e covariant net, the modular unitaries associated
with wedge regions should necessarily coincide with the boost
unitaries in the given representation  if the
modular covariance and the split property are
assumed. Here, we do not have any given unitary representation
of the Poincar\'e group from the start. Instead we recover a
canonical representation of the Poincar\'e group, which  is
unique if the split property holds, just by using the
modular  unitaries corresponding to all wedge algebras.

This step forward makes two significant changes. The first
concerns physics, since our result shows that the net itself
contains all the space-time symmetry informations, and that the
modular covariance property gives rise to a
canonical representation even when (non-split case) it is
not unique.

The second concerns mathematics, providing a new application  of
the group cohomology describing central extensions of groups.

 As is well known, central extensions of groups appear naturally
in physics, as shown for example by the Wigner theorems
[\ref(Simm)]  on the (anti-)unitary realization of the
symmetries in quantum physics.

 The class of results in this context which is relevant
for our purposes concerns the triviality of (suitable classes of)
central extensions of the universal covering of the Poincar\'e
group.
 One of the first results of this type is due to Michel
[\ref(Mich1)], but his hypotheses seems unapplicable to our
problem.

Another result in this direction is a particular case of the
theory of Moore [\ref(Moor1),\ref(Moor2)] on central
extensions of locally compact groups via polish groups.
The extensions described by Moore are continuous and open, and
his theorem on the universal extension concerns extensions by
closed subgroups of $\U(\H)$, the group of unitary operators on a
separable Hilbert space $\H$.

As a main step in our analysis, we study a class of extensions
with natural measurable properties, which we call weak Lie
central extensions. These extensions are given by a (not
necessarily closed) subgroup $G$ of $\U(\H)$, and we assume that
there are continuous one-parameter subgroups in $G$ which
correspond to a set of Lie-algebra generators (see section 2 for
the precise definition). We show that weak Lie central
extensions of a simply connected perfect Lie group are trivial,
and this applies to the universal covering of the Poincar\'e
group.

Our proof strongly relies on Moore theory. In fact, even though
the continuity properties of the central extensions described by
Moore are too restrictive for our purposes, we prove that his
theorem on the structure of the first and second order
(measurable) cohomology of a connected perfect Lie group
(see [\ref(Moor2)] or theorem 1.3) can be generalized to the weak
Lie central extensions (see theorem 1.7).
 The cases of the Poincar\'e group and  of the conformal
group find application in section 2.

Finally we mention that our results, together with theorem 2.3
in our previous paper [\ref(BGL)], give a complete
characterization of positive-energy conformal pre-cosheaves as
those where the modular groups associated with
double cones act geometrically as expected (see corollary 2.9).

The plan of the paper is the following: in the first section we
review some  known facts about cohomology and extensions of
groups, recall elements of Moore theory and  provide a
generalization of results needed in the  following.

 In the second section  we present our results  about the
generation of the unitary representation of the group of
space-time symmetries by modular groups, in case of local
algebras on a separable Hilbert space.

The last section contains an outlook.

 \section 1. Group extensions and cohomology
 \par
 We begin by shortly reviewing  basic elements of central
extensions of groups (see  [\ref(McLa1),\ref(Brow1)]). We first
deal with {\it algebraic central extensions}, i.e. with short
exact sequences of groups
 $$
1\to A\buildrel i\over\to G\buildrel\pi\over\to P\to
1\eqno(1.1)
 $$
 where $P$ is the group to extend and $A$ is a central
subgroup of $G$. The pair $(G,\pi)$ determines the
extension, and we often refer to it as a central
extension of $P$ (by $A = ker(\pi)$).
 \par
 Two extensions $(G_1,\pi_1)$, $(G_2,\pi_2)$ are called
{\it equivalent} if there is an isomorphism between $G_1$
and $G_2$ such that the following diagram commutes:
 $$
\matrix{
 &   & &        &     G_1     &        & &   & \cr
 &   & &\nearrow&             &\searrow& &   & \cr
1&\to&A&        &\updownarrow &        &P&\to&1\cr
 &   & &\searrow&             &\nearrow& &   & \cr
 &   & &        &     G_2     &        & &   & \cr}
\eqno(1.2)
 $$
 An extension of $P$ via $A$ is  {\it trivial} if it
is equivalent to the direct product $(A\times P,\pi)$ where
$\pi$ is the projection onto $P$, i.e. sequence (1.1) splits.

 A {\it section} of an extension $(G,\pi)$, i.e. a
map $s:P\to G$ such that $\pi\cdot s=id_P$, determines an
identification of $G$ with $A\times P$ as sets given by
 $g\to(s(\pi(g))^{-1}g,\pi(g))$. Then, the
multiplication rule on $A\times P$ is given by
 $$
(a,p)\cdot(b,q)=(ab\o(pq)^{-1},pq),\eqno(1.3)
 $$
where $\o(p,q) = s(q)s(pq)^{-1}s(p)$ satisfies the 2-{\it
cocycle} condition $\d_2\o=0$, where $\d_2$ is defined in
$(1.4)$.
 Conversely a 2-cocycle gives an associative multiplication on
$A\times P$, thus defining an extension.
 Two cocycles $\o_1$, $\o_2$ give rise to
equivalent extensions {\it iff} there is a group
automorphism of $A\times P$ which is compatible with
diagram (1.2) and intertwines the products given by $\o_1$
and $\o_2$. This means that there is a map $\f:P\to A$
such that the morphism between $(A\times P,\o_1)$ and
$(A\times P,\o_2)$ is given by $(a,p)\mapsto(a\f(p),p)$,
and $\o_1=\o_2\cdot\d_1\f$, with $\d_1$ defined in (1.4),
 i.e. $\o_2$ differs by $\o_1$ of a coboundary, showing
that equivalence classes of extensions are indeed
cohomology classes.
 \par
 The $n$-{\it cochains} $C^n(P,A)$ of $P$ with values in $A$
  are maps from $P^n$ to $A$, and the
$n$-th {\it coboundary map} $\d_n:C^n(P,A)\to
C^{n+1}(P,A)$ is given by
 $$
\eqalign{
\d_n f(p_1,\dots,p_{n+1})
&=f(p_2,\dots,p_{n+1})-f(p_1p_2,p_3,\dots,p_{n+1})\cr
&+f(p_1,p_2p_3,\dots,p_{n+1})+\dots\cr
&+(-1)^nf(p_1,\dots,p_np_{n+1})-(-1)^nf(p_1,\dots,p_n)}
\eqno(1.4)
 $$
and verifies $\d_{n+1}\d_n=1$.
 \par
 The range of $\d_{n-1}$ is denoted by
$B^n(P,A)$, the kernel of $\d_n$
by  $Z^n(P,A)$, and their quotient by
$H^n(P,A)$. Since $C^n(P,A)$ is an abelian group with respect to
the pointwise multiplication, $H^n(P,A)$ is a group too.
 \par
 The  above discussion shows that equivalence classes of
central extensions of $P$ via $A$ are in $1-1$
correspondence with elements of $H^2(P,A)$, and therefore
form a group. Note that the cocycle equation for
1-cochains means exactly that $Z^1(P,A)\equiv H^1(P,A)$
is the group of homomorphisms from $P$ to $A$, therefore
the vanishing of the first cohomology group for all
abelian $A$ is equivalent to the fact that $P$ coincides
with its commutator $[P,P]$, i.e. $P$ is {\it perfect}.
\par
 The groups $G$ for which all central extensions
split in a unique way, namely $H^n(P,A) =0,  n=1,2$ for all A as
above, are called {\it algebraically simply
connected} (see [\ref(Miln1),\ref(Brow1)]).
 \par
 A {\it universal central extension} $(E,\s)$ of $P$ is a
central extension such that for each central extension $(G,\pi)$
of $P$ there is a homomorphism $U:E\to G$ such that the
following diagram commutes:
 $$

\def\mapse#1{\searrow\hskip-.3cm {}^#1}
\def\mapright#1{\smash{
      \mathop{\longrightarrow}\limits^{#1}}}
\def\mapdown#1{\Big\downarrow
    \rlap{$\vcenter{\hbox{\hskip-.1cm $\scriptstyle#1\/$}}$}}
\matrix{
E &\mapright U&G\cr
&\mapse{\s}&\mapdown{\pi}\cr
&&P\cr}\ .\eqno(1.5)
 $$
 If a universal central extension exists, it is unique (up to
equivalence) and the extensions of $P$ via $A$ may be easily
described.
 Indeed consider the following commutative diagram:
 $$

\def\mapright#1{\smash{
     \mathop{\longrightarrow}\limits^{#1}}}
\def\mapdownu#1{\Big\downarrow
    \rlap{$\vcenter{\hbox{$\scriptstyle#1$}}$}}
\def\mapdownd#1{\Big\Vert
    \rlap{$\vcenter{\hbox{$\scriptstyle#1$}}$}}
\matrix{
&1&\mapright{}&A&\mapright{}&E\times
A&\mapright{id_E\cdot1_A}&E&\mapright{}&1&\cr
&&&\mapdownd{}&&\mapdownu{U\cdot id_A}&&\mapdownu\sigma
&&&\cr  &1&\mapright{}&A&\mapright{}&G&\mapright\p
&P&\mapright{}&1&.\cr}\eqno(1.6)
 $$
 By a classical diagram chasing argument, it is easy to
see that $U\cdot id_A$ is surjective, and that the
restriction of the projection $id_E\cdot1_A$ to
$\ker(U\cdot id_A)$ is an isomorphism onto $S:=\ker(\s)$.
Moreover the embedding of $S$ in $E\times A$ via the
preceding isomorphism has the form $id_S\times\psi$, where
$\psi$ is a morphism from $S$ to $A$.
 As a consequence, there is an isomorphism
 $$
\quot{E\times A}S \rightarrow G\eqno(1.7)
 $$
where $S$ is embedded into $E\times A$ as described,
therefore the extension $(G,\pi)$ is described in terms
of this embedding, i.e. in terms of the morphism
$\psi$, which implies that $H^1(S,A)\simeq H^2(P,A)$.
 \par
 A central extension $(E,\s)$ of $P$ where $E$ is
algebraically simply connected is the universal central
extension of $P$. Indeed fix an extension $(G,\pi)$
of $P$, choose a section $s$ of $\pi$ and consider the
2-cocycle $\o$ defined in equation (1.3). Then the map
 $$
 \matrix{
 \tilde\o&:E\times E&\to&A\cr
 &(g,h)&\mapsto&\o(\s(g)\s(h))}\eqno(1.8)
 $$
 is indeed a cocycle in $Z^2(E,A)$. Since $H^2(E,A)$ is
trivial, there is a 1-cochain $\f$ s.t.
$\d_1\f=\tilde\o$. Then, setting
 $$
U(g)=s(\s(g))\f(g)^{-1},\qquad g\in E\eqno(1.9)
 $$
 we easily get that $U$ is a homomorphism from
$E$ to $G$ and makes diagram (1.5) commutative.
 We note that since $\s$ is surjective and $E$ is
perfect, necessarily $P$ is perfect too.
\par
 If $P$ is a topological group,  extensions with given
topological properties are of interest. A theory in this sense
has been developed by Moore in a series of papers
[\ref(Moor1),\ref(Moor2)]. We recall that a topological group
$A$ is called {\it polish} if its topology may be obtained by
 a separable complete metric (that can be chosen compatible with
the uniformities of $A$).
  Moore considers {\it topological central extensions} of a
locally compact group $P$ via a polish group $A$, i.e. exact
sequences (1.1) such that $i$ is a homeomorphism into its image
and $\pi$ is continuous and open. He denotes by $Ext(P,A)$ the
equivalence classes of topological extensions of $P$ via $A$,
where the isomorphism in diagram (1.2) is asked to be a
homeomorphism. In analogy with the algebraic case, topological
cohomology groups are to be defined in such a way that the
identifications of $Ext(P,A)$ with $H^2_{top}(P,A)$, and of
$H^1_{top}(P,A)$ with the continuous homomorphisms from $P$ to
$A$ hold. Moore proves this to be the case if  cochains are
Borel measurable. With this hypothesis, the previously mentioned
identifications hold, cohomology groups are polish groups and
all usual functorial properties are satisfied. From now on we
consider topological cohomology groups only, therefore we drop
the subscript $_{top}$.
 \par
 Let us call {\it unitary group} a
subgroup of $\U(\H)$, the group of the unitary operators on a
separable Hilbert space $\H$ equipped with the weak topology.
It turns out that closed (i.e. complete) unitary groups  are
polish. Concerning universal extensions, Moore deals with {\it
unitary topological extensions}, namely $A$ in (1.1) is a closed
unitary group.
 We shall say that $E$ is $\U$-{\it simply connected} if all
central unitary extensions split in a unique way, that is to say
$H^1(E,A)$ and $H^2(E,A)$ vanish for all closed abelian unitary
groups $A$.
 \claim{1.1 Theorem [\ref(Moor2)]} Let the locally compact
group $P$ have a universal topological central extension $(E,\s)$
with $E$ a polish $\U$-simply connected group. Then
$(E,\s)$ is a universal central extension for the class of all
unitary central extensions, namely for
every topological central extension $(G,\pi)$ via a closed
unitary group $A$ there exists a continuous map $U:E\to G$ s.t.
the diagram (1.5) commutes. As a consequence, the isomorphism
(1.7) is a homeomorphism and $H^1(S,A)$ is isomorphic to
$H^2(P,A)$, where $S:=\ker(\s)$.
 \par
 According to theorem 1.1, we call the $\U$-covering of $P$ the
(unique) central extension of $P$ via a simply connected group
(if it exists). We need however a generalization of theorem 2 to
a larger class of extensions.
 \claim{1.2 Theorem}  Let $P$ be a locally compact group which
admits a $\U$-covering $(E,\s)$. Let $(G,\pi)$ be an algebraic
central extension of $P$ with $G$ a unitary group and
$A:=\ker(\pi)$ a topological subroup of $G$, closed as unitary
group. If $\pi$ has a Borel  section $s$, then $\pi$ splits on
$E$, i.e there exists a continuous map $U:E\to G$ s.t. diagram
(1.5) commutes.
 \np
 As a consequence, the group $\quot{E\times A}S$ is
polish and isomorphism (1.7) is continuous.
 \np
 If, in addition, the map $\pi$ is continuous, then $G$
is a closed unitary group and the extension $(G,\pi)$ is
a topological extension in the sense of definition 1.1.
 \par
 The last statement of the theorem shows that the
difference between the extensions in theorems 1.1 and 1.2 lies
in the continuity of the map $\pi$.
 \proof
 Since $s$ is a measurable section, the 2-cocycle defined
in (1.3) is measurable too, and it is an element of
$Z^2(P,A)$. Then we define $\tilde\o$ as in (1.8), and
since $\s$ is continuous, $\tilde\o$ turns out to be an
element of $Z^2(E,A)$. Since $A$ is a closed unitary group
and $E$ is $\U$-simply connected, $\tilde\o$ is indeed a
coboundary, i.e. there exists $\f\in C^1(E,A)$ s.t.
$\tilde\o=\d_1\f$. Then $U$ defined by equation (1.9)
is a homomorphism from $E$ to $G$, it is measurable by
construction, and makes diagram (1.5) commute. By
proposition $5(a)$ in [\ref(Moor1)], $U$ is indeed
continuous.
 \np
 Now we come back to diagram (1.6). Surjectivity of
$U\cdot id_A$ and the existence of the isomorphism
$j:\ker(U\cdot id_A)\to S$ follow by algebraic reasons.
Then we observe that $j$ is simply the restriction to
$\ker(U\cdot id_A)$ of the projection $id_E\cdot1_A$,
which is continuous and open, therefore $j$ is a
homeomorphism. Then $\quot{E\times A}S$ is a polish
group and the isomorphism (1.7) is continuous.
 \np
 If the map $\pi$ is continuous, then
 $$
g\in G\to\left(\s^{-1}\circ\pi(g)\ ,\
\left[U\circ\s^{-1}\circ\pi(g)\right]^{-1}\
g\right)\in E\times A
 $$
gives a map from the subsets of $G$ to the subsets of
$E\times A$ for which the preimage of an open set is
open. On the other hand, composing such a map with the
projection on $\quot{E\times A}S$, we get a map of
points which is the inverse of the isomorphism
(1.7). Such isomorphism is therefore a homeomorphism,
and all other properties follow.
 \endproof
 Note that, in the topological setting, the only
property that $[P,P]$ is dense in $P$ is needed for the
existence of a universal central extension. A remarkable
result of Moore shows that a perfect connected group $P$
admits the $\U$-covering.
Moreover, if $P$ is a Lie group, central extensions may be
completely described in terms of the Lie algebra cohomology and
the fundamental group $\pi_1(P)$.
 \par
 We recall that when $\P$ is a Lie algebra, a central
extension of $\P$ via a vector space $\V$ is an exact
sequence $0\to\V\to\G\buildrel\pi\over\to\P\to0$ where
$\V$ is a central Lie subalgebra of $\G$ and $\pi$ is
a Lie algebra homomorphism. Extensions are still
described in terms of 2-cocycles, and equivalence
classes of extensions in terms of 2-cohomology classes. We
denote by $H^n(\P,\V)$ the cohomology groups of $\P$
with respect to $\V$.
 \par
 \claim{1.3 Theorem [\ref(Moor2)]} The $\U$-covering
$(E,\s)$ of a perfect connected Lie group $P$ is a perfect
Lie group, and
 $$
\ker(\s)\simeq\pi_1(P)\times H^2(\P,\Re)
 $$
 where $\P$ is the Lie algebra of $P$.
 If $H^2(\P,\Re)=0$, $E$ coincides with the universal covering of
$P$.
 \par
 In the case of Lie groups, the existence of a measurable
section in theorem 1.2 may be replaced by a natural
condition.
 \claim{1.4 Definition} Let $(G,\pi)$ be an algebraic
extension of the Lie group $P$ where $G$ is a unitary
group. We say that $\pi$ is a weak Lie extension if there
exists a set $\{L_1\dots L_n\}$ of generators of $\P$ as
a Lie algebra ad a corresponding set
$\{V_1(t),\dots,V_n(t)\}$ of strongly continuous
1-parameter groups in $G$ such that
 $$
\pi(V_i(t))=\exp(tL_i),\qquad i=1,\dots,n.
 $$
where we have denoted by $\P$ the Lie algebra of $P$.
 \par
 \claim{1.5 Lemma} Let $(G,\pi)$ be a central weak Lie extension
of the Lie group $P$. Then there is a finite set of strongly
continuous 1-parameter groups in $G$ such that the generators of
their images in $P$ are a basis for $\P$ as a vector space.
 \par
 \proof Let us consider the map $v$ which associates to a
given set of one parameter groups in $P$ the (finite
dimensional) vector space spanned by their generators.
 \np
 We observe that given $A,B\in\P$, the commutator
$[A,B]$ belongs to
 $$
X:=v(\{t\to\exp(sA)\ \exp(tB)\ \exp(-sA):s\in\Re\}).
 $$
 Indeed, $X$ is spanned, by definition, by the range of
the Lie algebra valued function  $s\to\exp(sA)\ B\
\exp(-sA)$, whose derivative in 0 is $[A,B]$.
 \np
 As a consequence, by finite dimensionality, we may find
a finite set $\{s_i^j\}$ such that the vector space
 $$
v\circ\pi\left\{t\to V_{h_1}(s_1^j)\cdot\dots\cdot
V_{h_p}(s_p^j) \cdot V_{h_0}(t)\cdot V_{h_p}(s_p^j)
\cdot\dots\cdot V_{h_1}(s_1^j):j=1,2,\dots\right\}
 $$
contains the element
$[L_{h_1},[\dots[L_{h_p},L_{h_0}]\ ]$, where the $V_i$ are
described in definition 1.4. Then the thesis easily follows.
 \endproof
 \claim{1.6 Proposition} Let $(G,\pi)$ be a central weak
Lie extension of the connected Lie group $P$, with $G$ a
unitary group. Then there exists a Borel measurable
section $S:P\to G$ of the extension $\pi$.
 \par
 \proof According to the previous lemma, we may suppose
that $\{L_1,\dots L_n\}$ is a basis for $\P$ as a vector
space.
 \np
 The map $\a :\Re^n\to P$ given by
 $
\a (t_1,\dots t_n)=\prod_{i=1}^n \exp(t_iL_i)
 $
 has non trivial Jacobian at the origin because
$\{L_1\dots L_n\}$ is a basis. Therefore there there exist
two open sets, $U\subset\Re^n$ and
$V\equiv \a (U)\subset P$, where $U$ is a neighborhood
of zero in $\Re^n$ and, consequently, $V$ is a
neighborhood of the identity element in $P$,  such that
the map $\a \,:\, U\to V$ is a diffeomorphism.
 \np
 Now we define the (strongly) continuous map $\b:\Re^n\to
G$ given by $\b(t_1,\dots t_n)=\prod_{i=1}^n V_i(t_i)$,
and observe that the map $S_0:=\b\ \a^{-1}:V\to G$ is a
(strongly) continuous section of $\pi$ on $V$, i.e.
$\pi\cdot S_0\equiv id_V$.
 \np
 Since $V$ is an open neighborhood of the origin and $P$
is a connected group, then $P$ is algebraically generated
by $V$, i.e., each element of $P$ can be written as a
finite product of elements in $V$.
 \np
 Now we consider the open covering
 $\bigcup_{g\in P} gV= P.$
 Since $P$ is $\sigma$-compact,  we may extract a
countable sub-covering,
 $$
  P =\bigcup_{k\in\Na}g_k V.
 $$
 Finally, define the measurable
partition of $ P $ given by,
 $$
 \eqalign{ A_1 &= g_1 V\cr A_{n+1}&=\, g_{n+1}
V\cap\left(\bigcup_{k=1}^{n} g_k V\right )^c \quad
n\in\Na.\cr }
 $$
 \np
 For each $n\in\Na$, we may write $g_n =v_{n,1} \cdots
v_{n,m_n}$, $v_{n,k}\in V$, hence, since $A_n \subset g_n
V$, any element of $A_n$ may be written as $ v_{n,1}
\cdots v_{n,m_n} v$, when $v$ varies in $V$. Then, we
define
 $$
 \matrix{
S_n :&\quad A_n &\longrightarrow&G&\cr
&&&&\cr
&v_{n,1} \cdots v_{n,m_n} v &\longmapsto&S_0 (v_{n,1})
\cdots S_0 (v_{n,m_n})S_0 (v)&\quad.\cr
}
 $$
 As before, $S_n$ is a strongly continuous section of $\pi$
on $A_n$. Therefore it is easy to understand that the map
 $S=\sum_{n\in\Na}\chi_{A_n}\,\,S_n$, where $\chi_B$ is
the characteristic function of $B$, is the desired
measurable section.
 \endproof
 Making use of theorems 1.2 and 1.3 and proposition 1.6,
we can prove the following theorem on weak Lie
extensions. We remark that no condition on the closure of
$A$ is needed.
 \claim{1.7 Theorem} Let $P$ be a perfect connected Lie
group. Then all central weak Lie extensions $(G,\pi)$ of $P$,
with $G$ a unitary group, split on the $\U$-covering
$(E,\s)$, i.e. there exist a continuous homomorphism
$U:E\to G$ such that diagram (1.5) commutes.
 \par
 \proof Let $\tilde G$ be the group algebraically generated by
$G$ and $\ov A$, where $\ov A$ is the weak closure of $A$ in
$\U(H)$, and extend $\pi$ to a morphism $\tilde\pi$ in such a
way that $\tilde\pi|_{\ov A}=0$. Then,  applying  proposition~1.6
and theorem~1.2  to the extension $(\tilde G,\tilde\pi)$
we get a continuous map  $U:E\to\tilde G$ such that
$\tilde\pi\cdot U=\s$.
 \np
 By theorem 1.3 $E$ is perfect, therefore
 $U(E)=U([E,E])\subseteq[\tilde G,\tilde G]$.
 Since $\ov A$ commutes with $G$ and $\tilde G$,
equality $[\tilde G,\tilde G]\equiv[G,G]$ holds. As a
consequence, $Rg(U)\subseteq G$, and the proof is
concluded.
 \endproof
 We recall that the (4-dimensional) Poincar\'e group
$\Poi$ is perfect and its Lie algebra $\wp$ satisfies
$H^2(\wp,\Re)=0$. Therefore, next corollary immediately
follows.
 \claim{1.8 Corollary} Let $(G,\pi)$ be a central weak Lie
extension of the Poincar\'e group $\Poi$ where $G$ is a
unitary group.
 Then there exists a strongly continuous unitary
representation $U$ of the universal covering $\Spin$ of
$\Poi$ such that $\pi\cdot U=\s$, where $\s$ is the
covering map, and there is a continuous isomorphism
 $$
\quot{\Spin\times A}{\Ze_2}\to G
 $$
 where $\Ze_2$ is a suitable order two central subgroup of
$\Spin\times A$.
 \par
 \proof The corollary easily follows by theorems 1.2, 1.3 and
1.7, the previous observations on the Lie algebra $\wp$ and the
fact that the kernel of the covering map from $\Spin$ to $\Poi$
is $\Ze_2$.
 \endproof

\section 2. Modular covariance and the
reconstruction of space-time symmetries.
 \par
  In this section we study how the unitary representation of
the group of the space-time symmetries of a  Quantum Field
Theory [\ref(Haag1)] may be generated by the modular unitaries
associated with the algebras of a suitable class of regions.
Local Quantum Theories are described by a local pre-cosheaf of
von~Neumann algebras (see [\ref(GL1)]), i.e by a map
 $$
\A:\O\to\A(\O),\qquad \O\in\K
 $$
where $\K$ is the family of the double cones in the Minkowski
space $M$, such that
 $$
\eqalign{ \O_1\subset\O_2&\imply\A(\O_1)\subset\A(\O_2)
\hskip1.5cm {\rm(isotony)}\cr
\A(\O)&\subset\A(\O')'\hskip3.cm {\rm(locality)}.\cr}
 $$
 \par
 Local algebras are supposed to act on a {\it separable} Hilbert
space $\H$, and $\A_0$ denotes the $C^*$-algebra of quasi-local
observables generated by the local algebras. There is a vector
$\Q$ (vacuum) which is cyclic for all of them.
 \par
 The pre-cosheaf is extended by additivity to general open
regions of $M$.
 \par
As a consequence of locality, $\Q$ turns out to be cyclic
and separating for the algebras associated with all non
empty open regions whose complement has non-empty
interior (Reeh-Schlieder property).
 \par
 We recall that a {\it wedge} region is any Poincar\'e
transformed of the region $W_1:=\{x\in\Re^n:|x_0|<x_1\}$.
The {\it boosts} preserving $W_1$ are the elements of the
one-parameter subgroup $\L_{W_1}(t)$ of $\Poi$ which acts
on the coordinates $x_0$, $x_1$ via the matrix
 $$
\left(\matrix{\cosh2\pi t&-\sinh2\pi t\cr
-\sinh2\pi t&\cosh2\pi t}\right)
 $$
 and leaves the other coordinates unchanged.
 One-parameter boosts $\L_W(t)$ for any wedge $W$ are
defined by Poincar\'e conjugation. We denote by $\W$ the
family of all wedges in $M$.
 \par
 We make here our main assumption:
\claim{2.1 Assumption: Modular Covariance} Given any wedge region
$W$, the modular  unitaries $\D_W^{it}$ associated with
$(\A(W),\Q)$ act geometrically on the pre-cosheaf $\A$,
i.e.
 $$
\D^{it}_W\A(\O)\D^{-it}_W=\A(\L_W (t)\O),
\qquad\O\in\K,\quad W\in\W\ .
 $$
 \par
 \claim{2.2 Proposition} Let $\A$ be a local pre-cosheaf on
$M$ satisfying modular covariance. Then essential duality
holds, i.e.
 $$\A(W)'=\A(W')$$
for each wedge region $W$.
 \par
 \proof The proof is identical to that of Theorem
2.3\footnote{$^1$}{We take this opportunity to point out
that the positivity of the energy is implicitly assumed
in corollary 2.2 and theorem 2.3 of [\ref(BGL)] in order
to ensure the uniqueness of the conformal group
representation. We thank H. Borchers for a comment in
this sense. However, positivity of the energy is not used
in the derivation of the essential duality in [\ref(BGL),
theorem 2.3$(i)$]}, part~$(i)$ in  [\ref(BGL)].
 \endproof
 Now we consider the universal covering $\Spin$ of the
Poincar\'e group $\Poi$. Given a wedge $W$, we shall indicate
with $\tilde{\L}_W(t)$ the lifting of ${\L}_W(t)$ on $\Spin$. If
$g\to U(g)$ is a unitary representation of $\Spin$, we shall
say that it is $\A$-covariant if
 $$
 U(g)\A(\O)U(g)^*=\A(\s(g)\O),\qquad g\in\Spin
 $$
 where $\s:\Spin\to\Poi$ is the covering map, and that it
has positive energy if the restriction of $U$ to the
translation subgroup satisfies the spectrum condition.
 \claim{2.3 Theorem} Let $\A$ be a modular covariant, local
pre-cosheaf on $M$, dim$(M)>2$.
Then there exists a vacuum preserving, positive energy,
$\A$-covariant, unitary representation $U$ of $\Spin$ which is
canonically determined by the property
 $$U(\tilde\L_W (t))=\D_W^{it},\qquad W\in\W$$
 \par
 In order to prove Theorem $2.3$ we need some results
about the structure of $\Poi$ (Lemma $2.5$), and about
the group generated by modular operators associated with
wedge regions.
 \claim{2.4 Proposition}  The generators of the boost
transformations generate the Lie algebra of the Poincar\'e
group. As a consequence, the boosts algebraically generate
$\Poi$.
 \par
 \proof The first statement is well known, and the second
follows by the connectedness of $\Poi$.
 \endproof
 \claim{2.5 Lemma} Let $\A$ be a local pre-cosheaf on
$M$ satisfying assumption $2.1$, and consider the set
 $$
H=\{g\in\Poi:\A(g\O)=\A(\O)\quad\forall\O\in\K\}
 $$
 Then, H is a normal subgroup of $\Poi$ and either
\item{$(a)$} $H=\{1\}$
\np
or
\item{$(b)$} for each $\O\in\K$, $\A(\O)$ is equal to
$\A_0$ and is a maximal abelian subalgebra of $B(\H)$.
 \par
 \proof $H$ is clearly a group, and we prove that it is indeed
normal in $\Poi$. In fact, by modular covariance, $\forall
g\in H$ we  get
 $$
\A(\L_W (t) g\L_W (-t)\O)=Ad\D_W^{it}\cdot\A(g\L_W (-t)\O)=
Ad\D_W^{it}\cdot\A(\L_W (-t)\O)=
\A(\O)
 $$
 therefore, since the boosts generate $\Poi$
(proposition~2.4), $H$ is normal.
 \np
Now we prove that each non trivial normal subgroup $K$ of
$\Poi$ contains the translation subgroup
$T=\{T(a):a\in\Re^n\}$.  Let $1\not=g\in K$. If $g$ is not
a translation then, by normality,
 $$
K\ni g T(a) g^{-1} T(-a)=T(g(a))T(-a)\in T
 $$
i.e. $K\cap T$ contains at least a non trivial element
$g$. Again by normality, $K\cap T$ is globally invariant
under the action of the Lorentz group. Since the (minimal)
Lorentz-invariant subsets of $\Re^n$ are the future light
cone, its surface, the past light cone, its surface, and
the complement of the preceding ones, and since any of
these generates $\Re^n$ as a group, then $K\supset T$. In
particular, we proved that if $H$ is non trivial it
contains $T$.
 In this case, for each $\O\in\K$ we can find a
translation $T(a)$ such that $T(a)\O$ is spatially
separated by $\O$, and therefore, by locality, $\A(\O)$ is
abelian. Then it is easy to see that all local algebras
coincide, and they are maximal abelian sub-algebras
[\ref(StraZs1)]  because the vacuum is cy\-clic.
 \endproof
 \rmclaim{Remark} If alternative $(b)$ of Lemma 2.5 holds, then
Theorem 2.3 is true by taking the trivial representation of the
universal covering of the Poincar\'e group.
 In the rest of the paper we shall discuss alternative $(a)$,
i.e. we shall suppose $H=\{1\}$.
 \par
 \claim{2.6 Lemma} Let $G$ be the subgroup of $\U(\H)$ which is
algebraically generated by the modular groups of the algebras
$\A(W)$, $W\in\W$. Then $G$ is an algebraic central extension of
the Poincar\'e group $\Poi$.
 \par \proof First we prove that the map $\pi$ which satisfies
 $$
 \p (\D^{it}_{W})=\L_{W}(t)\qquad W\in\W,\quad t\in\Re\eqno(2.1)
 $$
 extends to a well defined surjective group homomorphism
$\p :G\to\Poi$.
 \np
 Indeed let
$\D_1^{it_1}\cdot\dots\cdot \D_n^{it_n}=1$ be a non
trivial identity in $G$, where
 $\D_i$ is the  modular operator of $\A(W_i)$ and
$\{W_i , i=1,\dots,n\}$ is any collection of $n$ wedges
in $\cal W$. Then
 $$
\A(\O)=
\D_1^{it_1}\cdot\dots\cdot\D_n^{it_n}\A(\O)
\D_n^{-it_n}\cdot\dots\cdot\D_1^{-it_1}=
\A(\p (\D_1^{it_1}\cdot\dots\cdot\D_n^{it_n})\O),
\quad\forall\O\in\K.
 $$
 From the Lemma 2.5 we get
$\p(\D_1^{it_1}\cdot\dots\cdot \D_n^{it_n})=1$.
 As a consequence $\p$ is a well-defined homomorphism.
 Since the boosts algebraically generate $\Poi$
(proposition 2.4), $\pi$ is surjective.
 \np
 Now we observe that if $U\in\ker(\p)$, then
$\ad U(\A(W))=\A(\p(U)W)=\A(W)$ for each wedge $W$, therefore
$U$ commutes with the modular group of any wedge algebra
[\ref(StraZs1)], and hence with any element in $G$. As a
consequence, the exact sequence  $1\lar\ker(\p)\lar
G\lar\Poi\lar 1$ gives the announced central extension.
 \endproof
 \proofof{Theorem 2.3} Equation (2.1) and proposition
2.4 implies that the extension $\pi$ described in
lemma~2.6 is a weak Lie extension. Since dim$(M)>2$, the
Poincar\'e group is perfect and its Lie algebra has
trivial second cohomology. Therefore, by lemma 2.6 and
corollary 1.8, we get a strongly continuous unitary
representation $U$ of $\Spin$. Setting
$U_W(t):=U(\tilde\L_W(t))$, $W\in\W$, $t\in\Re$, the
unitary operator
 $$
z(t)=\D_W^{it}U_W(-t)
 $$
 implements internal symmetries by modular covariance,
therefore is in the center of $G$.
 Let us denote by $\th_g$, $g\in\Poi$, the action of
$\Poi$ on the central extension $G$.
It follows that $\th_g(\D_W^{it})=\D_{gW}^{it}$,
$\th_g(U_W(t))=U_{gW}(t)$. Since $\th_g$ acts trivially on
the center of $G$ and $\Poi$ acts transitively on the
family $\W$, then $z(t)$ does not depend on $W\in\W$.
Moreover
 $$
z(s)z(t)=\D_W^{it}z(s)U_W(-t)=z(s+t)\qquad s,t\in\Re
 $$
i.e. $z(t)$ is a one-parameter group. Finally, since
$\D_{W'}^{it}=\D_W^{-it}$ by essential duality (theorem
2.2), and $U_{W'}(t)=U_W(-t)$ because $\tilde\L_W (t)$ and
$\tilde\L_{W'}({-t})$ are the unique lifting of  $\L_W
(t)=\L_{W'}({-t})$,  we have $z(t)=z(-t)$, i.e.
$z(t)\equiv I$.
 \np
In order to check the positivity of the energy-momentum, we
observe that the generator of any time-like translation
is a convex combination of generators of light-like
translations. Moreover  each generator of a light-like
translation gives rise to a one-parameter  semigroup of
endomorphisms of a suitable wedge. Therefore the result follows
by proposition 2.7.
 \endproof
 The following proposition is a partial
converse to Borchers theorem, in a slightly different form it
appeared also in Wiesbrock paper [\ref(Wies1)], with the
difference that we do not assume the commutation relations of the
modular conjugation with the one-parameter group.
 \claim{2.7 Proposition} (Converse of Borchers theorem) Let $\R$
be a von Neumann algebra standard with respect to the vector
$\Q$. If $U(a)$ is a one-parameter semigroup ($a\in \Re_+ $)  of
endomorphysms of $\R$ such that
 \item{$(i)$}  $U(a)\Q=\Q$
 \item{$(ii)$}  $ad\D_{\R}^{it} U(a)=U(e^{-2\pi t} a)$
 \np
 then $U(a)$ has positive generator.
 \par
 We give here an easy proof of the preceding Theorem.
 \proof From Lemma II.$3$ in [\ref(Borc1)] we know that
the operator valued function $t\to U(e^{-2\pi t} a)$ can
be analitically extended on the strip $\{ z\in\Co ;
-1/2<Imz<0\}$, where it satisfies the bound $\|U(e^{-2\pi
z}a)\|\le 1$. Then, taking $z=-i/4$, we get
 $$
\|U(ia)\|=\|e^{-aP}\|\le 1
 $$
 hence $P$, being the generator of $U$, has non-negative
spectrum.
 \endproof
 \claim{2.8 Corollary} If the distal split property (cf.
[\ref(DoLo1)]) holds, the representation described in
Theorem $2.3$ is the only $\A$-covariant unitary
representation of the universal covering group $\Spin$ of
$\Poi$ on $\H$.
 \par
 \proof The proof is identical to that of Theorem $3.1$
in [\ref(BGL)].
 \endproof
 \par
 We conclude this section observing that our techniques
work also for conformally covariant
precosheaves (see e.g. [\ref(BGL)] for definitions and
properties). Moreover, in this case, the (conformal) modular
covariance property for a local precosheaf on (some covering of)
the Dirac-Weyl compactification of the Minkowski space is a
necessary and sufficient condition for the existence of a
covariant unitary representation of the conformal group.
 \par
 For simplicity, we illustrate this for conformal theories
on $S^1$. Such theories are described by a local precosheaf of
von~Neumann algebras on the family $\K$ of all proper open
intervals in $S^1$. Here the causal complement $I'$ denotes the
interior of the complement of $I$. As before, the vacuum vector
is supposed to be cyclic for the local algebras associated with
subsets in $\K$.
 \par
 With each $I\in\K$ one associates the one-parameter
subgroup of the M\"obius group $\L_I(t)$ which preserves
$I$ [\ref(BGL)]. The following theorem holds:
 \claim{2.9 Theorem} Let $I\to\A(I)$ be a local precosheaf of
von~Neumann algebras  on $S^1$ and $\Q$ a cyclic vector for
each $\A(I)$. Then $\A$ is conformally covariant if and only if
the modular groups associated to $(\A(I)$, $\Q)$, $I\in\K$,
verify (conformal modular covariance)
 $$
\D_I^{it}\A(L)\D_I^{-it}=\A(\L_I(t)L), I,L\in\K.
 $$
In this case the unitary representation $g\to U(g)$ of the
M\"obius group $SL(2,\Re)/\Ze_2$ is generated by the modular
groups of the algebras associated with open intervals and
$U(\L_I(t))=\D_I^{it}$, $I\in\K$.
 \par
 \proof With the same arguments used in the Poincar\'e
covariant case, the unitary group generated by the modular
automorphisms associated with double cones generates a weak Lie
central extension of the M\"obius group which, because of
theorem 1.7, splits on the universal covering. Again, we get a
unitary representation of the covering group of the M\"obius
group where the action of  one-parameter group which preserves an
open interval is implemented by the corresponding modular group.
Then, applying the same argument as in [\ref(BGL),\ref(FrGa1)],
the modular conjugations associated with open intervals implement
the reflection with respect to the end points of the interval.
An argument of Wiesbrock [\ref(Wies2)] finally implies that the
constructed representation is indeed a representation of the
M\"obius group.
 \endproof
 \section 3. Final comments.
 \par
 In a forthcoming paper [\ref(FGL)]
 our analysis will be completed in two directions. First,  a
discussion on the relationship between the modular covariance
(assumption 2.1) and the geometrical behaviour of the modular
conjugations associated with wedge regions will be discussed,
with applications to the construction of a PCT symmetry.

 Secondly, notice that we have obtained a representation of the
universal covering group $\Spin$ of the Poincar\'e group $\Poi$,
therefore the center $\Ze_2$ of $\Spin$ might act as a
non--trivial gauge symmetry. Our net being local, this
possibility is incompatible with the spin-statistics
correspondence [\ref(Haag1)] and will be analyzed in terms of an
algebraic spin-statistics theorem.  Moreover, we expect the
latter to follow from the split property [\ref(DoLo1)], that
selects physically relevant nets and implies the uniqueness of
the Poincar\'e action.

 Note however that the product of the modular conjugations
associated with an inclusion of wedge regions is necessarily a
translation [\ref(Borc1)], and our analysis should be
compared with the characterization of translation invariant
theories given in [\ref(BuSu1)].

 Finally we mention that
essential duality may fail in general, as shown recently in
[\ref(Yngv1)].

  \section Acknowledgements
\par
We would like to thank Prof. Sutherland for pointing out
the relevance of references [\ref(Moor1),\ref(Moor2)] for our
purposes. One of us (R.L.) is also grateful to him for the
warm hospitality at the University of New South Wales, Sidney,
during January-February 1993. The first named author is grateful
to Prof. Balachandran for the kind invitation at the Department
of Physics of Syracuse University and to Prof. G. Marmo for
pointing out reference [\ref(Mich1)].

\references
\end